\documentclass[12pt]{JHEP3}
\usepackage{mathrsfs}
\usepackage{amsmath}
\usepackage{epsfig}
\title{ The generating function of the $\sigma_1$ function}
\author{Jingbo Wang \\
  Department of Applied Physics,
 Xi'an Jiaotong University,
 Xi'an, 710049, People's Republic of China\\
 \email{ shuijing31@gmail.com }}

\abstract { In this paper, we give a generate function for the
$\sigma_1$ function. Then we find some connections between the
$\sigma_1$ function and the Ramanujan's tau function. We hope this
connection will give some insights into the unsolved problems in
classical number theory. }

 \preprint{} \dedicated{}
\maketitle
\begin{document}
\section{Introduction}
Consider the following function:
\begin{equation}\label{1}
   f(x) = \sum\limits_{n = 1}^\infty  {\ln (1 - x^n )},\hskip 5mm
   0<x<1.
\end{equation}
With the help of the Taylor's expansion, we get:
\begin{equation}\label{2}
    f(x) = \sum\limits_{n = 1}^\infty  {\ln (1 - x^n )}  =  -
\sum\limits_{n = 1}^\infty  {\sum\limits_{m = 1}^\infty
{\frac{{x^{nm} }} {m}} }  =  - \sum\limits_{n = 1}^\infty
{\frac{{\sigma _1 (n)}} {n}} x^n.
\end{equation}
where the arithmetic function $\sigma_n$ is defined as follows(see
\cite{Na1}):
\begin{equation}\label{3}
    \sigma _n (m) = \sum\limits_{d|m} {d^n }.
\end{equation}
In this paper, we just consider $\sigma_1$, that is the sum of
positive divisors of $m$, and will omit the subscript $1$. Denote
$E(n)=\sigma(n)/n$, then $E(n)$ has some simple properties:
\begin{itemize}
\item $E(n)>1$ for every n;
\item when $p$ is a prime, $E(p)=(p+1)/p$;
\item there exist $n$ such that $E(n)=2$, and those numbers are called prefect
numbers;
\item $E(n)$ has no up bound. For example, consider $N=n!$, then
$E(N) > 1 + \sum\limits_{n = 2}^n {1/n}.$ When $ n \to \infty$,
$E(N) \to \infty$.
\end{itemize}
From the equation\ref{2} we can get that
\begin{equation}\label{4}
    - f(x) =  - \sum\limits_{n = 1}^\infty  {\ln (1 - x^n )}  = \sum\limits_{n = 1}^\infty  {E(n)}
    x^n.
\end{equation}
That is, $-f(x)$ is the generating function of the $E(n)$ function.
\section{The properties of the function $f(x)$}
Since $ln(1-x^n)<0$,
\begin{equation}\label{5}
 f(x) = \sum\limits_{n = 1}^\infty  {\ln (1 - x^n )}  =
\sum\limits_{n = 1}^\infty  {\ln (1 - \exp (n\ln x))}  > \frac{1}
{{\ln x}}\int\limits_{ - \infty }^0 {\ln (1 - \exp (y))dy}  =
\frac{{\zeta (2)}} {{\ln x}}.
\end{equation}
On the other hand, $E(n)>1$, we get
\begin{equation}\label{6}
   f(x) =  - \sum\limits_{n = 1}^\infty  {E(n)} x^n  <  -
\sum\limits_{n = 1}^\infty  {x^n }  = \frac{x} {{x - 1}}.
\end{equation}
Combine those two equation, we get
\begin{equation}\label{7}
\frac{{\zeta (2)}} {{\ln x}} < f(x) =  - \sum\limits_{n = 1}^\infty
{E(n)} x^n  < \frac{x} {{x - 1}}, \hskip 5mm 0<x<1.
\end{equation}
In the following graph, we plot those three functions. Since the
function $f(x)$ has infinity terms, we just add the first ten and
twenty terms, so it intersect with function $x/(x-1)$. If we add
more and more terms, the $f(x)$ will be between the other two
function better and better. \EPSFIGURE{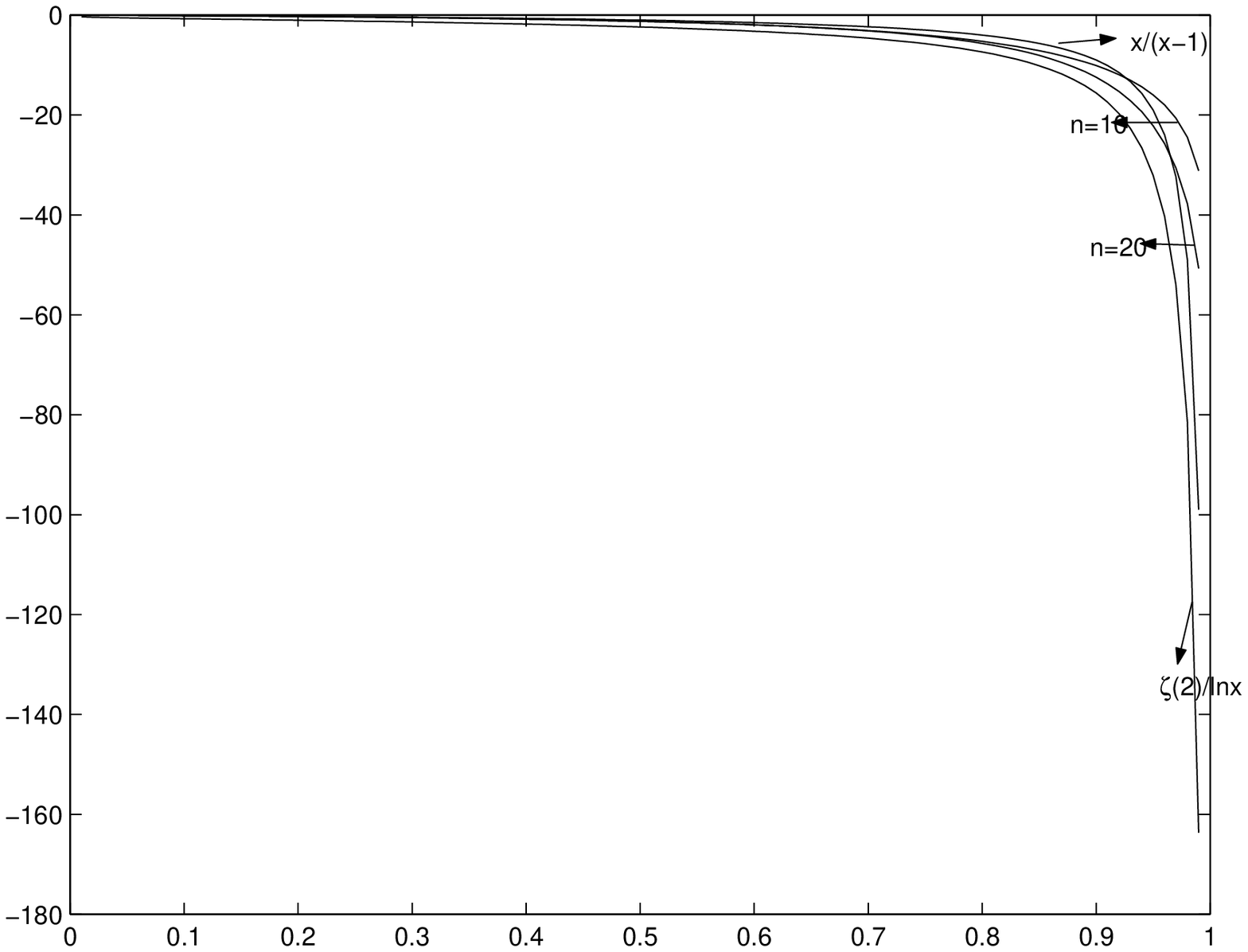} {Three
functions in $0<x<1$}

The right side of the inequality seems trivial, but the left side
may give some constrains on the arithmetic function $E(n)$ and the
related $\sigma(n)$.
\section{Relation with Ramanujan's tau function}
The Ramanujan's tau function is defined implicitly by\cite{Se1}
\begin{equation}\label{8}
   x\prod\limits_{n = 1}^\infty  {(1 - x^n )^{24} }  = \sum\limits_{n =
1}^\infty  {\tau (n)x^n }.
\end{equation}
From the above equation we can get the relation between the $E(n)$
and the $\tau$.
\begin{equation}\label{9}
   \sum\limits_{n = 1}^\infty  {\tau (n)x^n }  = x\prod\limits_{n =
1}^\infty  {(1 - x^n )^{24} }  = \exp (\ln x + 24\sum\limits_{n =
1}^\infty  {\ln (1 - x^n )} ) = x\exp ( - 24\sum\limits_{n =
1}^\infty  {E(n)} x^n ).
\end{equation}
Expanding the exponential function on the right side, we can get the
formulas relate the $E(n)$ and the $\tau$. For example, we can get
$\tau(1)=1, \tau(2)=-24E(1)=-24,
\tau(3)=-24E(2)+1/2*24^2*E(1)=252,\cdots $just the right numbers.
Unfortunately we can't get the general formula to calculate the
$\tau(n)$ from $E(n)$, or vice verse. From the Ramanujan's tau
function we can get the simplest cups form. They also satisfy the
Ramanujan's conjectures(established by Delinge), that is, $\left|
{\tau (p)} \right| \leqslant 2p^{11/2}$ for all primes p. This can
also give some constrains to the function $E(n)$.
\section{Conclusion}
In classical number theory, there are many unsolved
problems\cite{Na1}. The related problems for this paper contain "are
there infinite many even prefect numbers?" and "are there any odd
prefect numbers?" and so on. In this paper, we relate those problems
to modern arithmetic, such as the modular forms, L-function and so
on. We hope those relations can give some insight into those
problems. \acknowledgments{This work was partly done at Beijing
Normal University. This research was supported in part by the
Project of Knowledge Innovation Program (PKIP) of Chinese Academy of
Sciences, Grant No. KJCX2.YW.W10 }
\bibliography{number1}
\end{document}